\documentclass[12pt]{amsart}
\usepackage{amssymb,latexsym}

\numberwithin{equation}{section}
\newtheorem{thm}{Theorem}[section]

\newtheorem{prop}[thm]{Proposition}
\newtheorem{cor}[thm]{Corollary}
\newtheorem{rem}[thm]{Remark}
\newtheorem{lem}[thm]{Lemma}

\newtheorem{definition}[thm]{Definition}

\usepackage{graphics, graphicx}

\begin{document}

\title[KOSZULITY OF SPLITTING ALGEBRAS]
{KOSZULITY OF SPLITTING ALGEBRAS ASSOCIATED WITH CELL COMPLEXES}

\author{Vladimir Retakh}
\address{Department of Mathematics, Rutgers University, Piscataway, NJ 08854-8019, USA}
\email{vretakh@math.rutgers.edu}
\author{Shirlei Serconek}
\address{IME-UFG, CX Postal 131, Goiania - GO, CEP 74001-970,  Brazil}
\email{serconek@math.rutgers.edu}
\author{Robert Lee Wilson}
\address{Department of Mathematics, Rutgers University, Piscataway, NJ 08854-8019, USA}
\email{rwilson@math.rutgers.edu}

\keywords{Splitting algebras, Koszulity, cell decompositions}
\subjclass{16S37; 57M20; 57M60; 55U10}

\begin{abstract}
We associate to a good cell decomposition of a manifold $M$ 
a quadratic algebra and show that the Koszulity of the algebra 
implies a restriction on the Euler characteristic of $M$.  
For a two-dimensional manifold $M$ the algebra is Koszul 
if and only if the Euler characteristic of $M$ is two.

\end{abstract}
\maketitle

\section{INTRODUCTION}
\medskip
Let $\Gamma = (V,E)$ be a layered graph (where $V$ is the set of 
vertices and $E$ is the set of edges).  One may define an assoicated 
algebra, $A(\Gamma)$, to be the algebra generated by $E$ 
subject to the relations which state that
$$(t-e_1)(t-e_2)\dots (t-e_k)=
(t-f_1)(t-f_2)\dots (t-f_k)$$
whenever the sequences of edges $e_1,e_2,...,e_k$ and  $f_1,f_2,...,f_k$ 
define directed paths with the same origin and the same end
and $t$ is an independent central variable.  As this algebra
records information about the factorization of polynomials associated to
$\Gamma$ into linear factors, we call $A(\Gamma)$ the {\it splitting algebra} associated to $\Gamma$.

These algebras have been introduced and studied in \cite {GRW, GGRW, RSW1, RSW2, RSW3, RW}.

The algebra $A(\Gamma )$ is defined by a set of homogeneous relations.
For certain graphs $\Gamma$ (i.e, uniform graphs as defined in \cite{RSW1})
these relations are consequences of a family of quadratic relations and so the splitting
algebra $A(\Gamma)$  possesses a 
quadratic dual algebra $A(\Gamma )^!$. In \cite{RSW1} we  claimed that any algebra
$A(\Gamma )$ defined by a uniform layered graph is Koszul. Later, T.
Cassidy and B. Shelton constructed a counter-example to this statement and
this example forced us to rethink the situation. As a result, we
discovered a large class of non-Koszul algebras $A(\Gamma )$ where $\Gamma
$ is a graph defined by a regular cell complex $C$. Indeed, in Sections 4 and 5 we show
that, if $\Gamma $ is the layered graph defined by a ``good" cell
decomposition of a two-dimenional manifold $M$, then $A(\Gamma )$ is 
Koszul if and only if the Euler
characteristic of $M$ is two.

The analysis of \cite{RSW1} 
does provide a sufficient condition for $A(\Gamma)$ to be Koszul.  
We discuss this in Section 6 and use this sufficient condition to show 
that if $\Gamma$ is either a complete layered graph or the graph associated to an 
abstract simplicial complex then $A(\Gamma)$ is Koszul.

Actually, any cell subdivision $D$ of a manifold $M$ defines two
graphs:
$\Gamma_D $ with vertices $\sigma \in D \cup \{\emptyset,M\}$ and
the subgraph $\Gamma _D'$ with vertices $\sigma \in D \cup \{\emptyset\}$.

The properties of $A(\Gamma_D )$ and $A(\Gamma_D ')$ are quite different.
We prove that $A(\Gamma _D')$ is Koszul for any
simplicial complex $D$ but as mentioned before the Koszulity of
$A(\Gamma _D)$ imposes a strong topological condition on $M$.

We also consider a condition which is weaker than Koszulity.   
We call a quadratic algebra $A$ numerically Koszul if
the Hilbert series of $H(A,z)$ and $H(A^!,z)$ are related by the famous
formula $H(A,z)H(A^!,-z)=1$.
This is not a common terminology: sometimes algebras with similar 
properties are called
quasi-Koszul (see \cite{Pi}). Koszul algebras are numerically Koszul but
the converse is not true (see \cite{Ro, Po}). We prove (Section 3), however, that if 
the height of $\Gamma $ is less or equal to $4$ then the numerical 
Koszulity of $A(\Gamma )$ implies its Koszulity. We use this in proving the results of Section 5.

In our discussions of numerical Koszulity we need efficient 
techniques for computing the Hilbert series of $A(\Gamma)$ and $A(\Gamma)^!$.
We give such a technique (Lemma 1.3) for $A(\Gamma)$ in Section 1 and use this result to show that
the example of Cassidy and Shelton is not numerically Koszul.  
Section 2 is devoted to results about the Hilbert series of $A(\Gamma)^!$

We are grateful to T. Cassidy and B. Shelton for providing us with their counter-example.
  
\section{SPLITTING ALGEBRAS ASSOCIATED WITH LAYERED GRAPHS}

Recall the definition of a {\it layered graph}.
Let  $\Gamma = (V, E)$ be
a directed graph. That is, $V$ is a set (of vertices), $E$
is a set (of edges), and ${\bold t}: E \rightarrow V$ and ${\bold
h}: E \rightarrow V$ are functions. ${\bold t}(e)$ is the {\it
tail} of $e$ and ${\bold h}(e)$ is the {\it head} of $e$.

A vertex $u\in V$ is called {\it maximal } if there is no $e\in E$
such that ${\bold h}(e)=u$. A vertex $v\in V$ is called {\it
minimal } if there is no $e\in E$ such that ${\bold t}(e)=v$.

We say that $\Gamma$ is {\it layered} if $V$ is the disjoint union
of $V_i$, $0\leq i\leq n$, $E$ is the disjoint union of $E_i$, 
$1\leq i\leq n$,
 ${\bold t}: E_i \rightarrow V_i$, \ ${\bold h}: E_i \rightarrow V_{i-1}$.
We will write $|v|=i$ if  $v\in V_i$. In this case the number $i$
is called the {\it level} or the {\it rank} of $v$. Note that a layered 
graph has no loops. There is an obvious bijection between layered graphs 
and ranked partially ordered sets.

If $v, \ w \in V$, a path of the length $k$ from $v$ to $w$ is a sequence 
of
edges $\pi = \{ e_1, e_2, ...,e_k \}$ with ${\bold t}(e_{1}) =
v$, ${\bold h}(e_k) = w$ and ${\bold t}(e_{i+1}) = {\bold
h}(e_i)$ for $1 \leq i < k$.  We write $v = {\bold t}(\pi)$, $w =
{\bold h}(\pi)$ and call $v$ the tail of the path and $w$ the
head of the path. We also write $v > w$ if there is a path from
$v$ to $w$.

Recall the definition of a uniform layered graph.
Let $\Gamma$ be a layered graph. For $v \in V$ define $S_-(v)$ to be the set 
of all vertices  $w \in V$ covered by $v$,
i.e. such that there exists an edge with the tail $v$ and the
head $w$. 
For $v \in V_j, j \ge 2$, let $\sim$ denote the equivalence
relation on $S_-(v)$ generated by $u \sim w $ if 
$S_-(u) \cap S_-(w) \ne \emptyset.$
The following definition was introduced in \cite{RSW1}.
 
\begin{definition} The layered graph $\Gamma =(V,E)$ is said to be
uniform if,
for every $v \in V_j, j \ge 2$,
all elements  of $S_-(v)$ are equivalent under $\sim$.
\end{definition}

In \cite{GRSW} for each layered graph $\Gamma =(V,E)$ we constructed an 
associative algebra generated by the edges of $\Gamma $. 
Let $T(E)$ denote the free associative algebra  on $E$
over a field $K$. We are going to introduce a quotient algebra of 
$T(E)$.
We will do this by equating coefficients of polynomials associated
with pairs of paths with the same origin and the same end.
For a path $\pi = \{ e_1, e_2, ...,e_m \}$ define
$$P_{\pi}(t) = (t-e_1)...(t-e_m) \in T(E)[t] .$$

Note that $P_{\pi_1\pi_2}(t) = P_{\pi_1}(t)P_{
\pi_2}(t)$ if ${\bold h}(\pi_1) = {\bold
t}(\pi_2)$. Write
$$P_{\pi}(t) = \sum_{j=0}^{m} (-1)^{m-j} e(\pi,j)t ^j .$$

\begin{definition} Let
$R_E$ be the linear space in $T(E)$ spanned by
$$\{ e(\pi_1,k) - e(\pi_2,k) \ | \ k\geq 1,
{\bold t}(\pi_1)={\bold t}(\pi_2),
{\bold h}(\pi_1)={\bold h}(\pi_2)\}.$$
\end{definition}

Denote by $<R_E>$ the ideal generated by $R_E$ and set $$A(\Gamma) = T(E)/<R_E>.$$

Suppose that a layered graph $\Gamma =(V,E)$ of the height $n$ has a 
unique maximal 
vertex $M\in V_n$ and a unique minimal vertex $*\in V_0$.
To any path $\pi$
from $M$ to $*$ there corresponds a monic polynomial $P_{\pi}$ of degree $n$.  
Then the image $P_{\Gamma}(t)\in A(\Gamma )[t]$ of the polynomial 
$P_{\pi}(t)$ does not depend on the choice of $\pi$ and
any path from $M$ to $*$ corresponds a factorization of 
$P_{\Gamma}(t)$ into a product of linear terms. 

Therefore, we call $A(\Gamma )$ a {\it splitting algebra}
of $P_{\Gamma}(t)$. There is an injection of 
the group of level preserving automorphisms of 
$\Gamma $ into the  group of automorphisms  of $A(\Gamma)$).

From now on we suppose that $\Gamma =(V,E)$ has a unique minimal vertex 
$*$, i.e. $V_0=\{*\}$.  
It was shown in \cite{RSW1} that if $\Gamma$ is uniform then $A(\Gamma)$ 
is quadratic. In fact, this algebra $A(\Gamma )$ can be defined as the
algebra generated by the linear space $KV$ subject to certain relations 
$R_V\subset V\otimes V$ (see \cite{8}). So 
the dual algebra $A(\Gamma)^!$ is defined as a quadratic
algebra with the space of generators $V'$ and the space of relations
$R_V^{\perp}\subset V'\otimes V'$ where $V'$ is dual to $V$ and
$R_V^{\perp}$ is the annihilator of $R_V$. 

Recall that for any graded algebra $B=\oplus _{n\geq 0}B_n$ with finite 
dimensional $B_n$ one can define its Hilbert series (graded dimension)
as $H(B,z)=\sum dim\ B_nz^n$. Therefore, Hilbert series can be defined for 
a quadratic algebra $A = T(V)/<R_V>$ and for its dual algebra $A^!$. The Hilbert 
series for splitting algebras associated with layered graphs were 
constructed in \cite{RSW3} (see \cite{RW} for a more general result).

Let $\Gamma$ be a layered graph.  
We will write $H(A(\Gamma),z)$ in a form which is convenient for computation. 
Let 
$$s_{a,b} = \sum_{v_1 > ... > v_l > *, |v_1| = a,|v_l| = b} (-1)^l.$$

\begin{lem} $H(A(\Gamma),z)^{-1} =  1 + \sum_{i\ge 0}(\sum_{a \ge i \ge a-b+1} s_{a,b})z^i.$
\end{lem}

\begin{proof}

By \cite {GRW} we have
$$(H(A({\Gamma},z))^{-1} = \frac {1 -z + \sum_{v_1 > ... > v_l \ge *} (-1)^l(z^{|v_1| - |v_l| + 1} - z^{|v_1| + 1})} {1-z}.$$
As $|*| = 0$, this is equal to
$$ 1+  \sum_{v_1 > ... > v_l > *} (-1)^l z^{|v_1| - |v_l| + 1} \frac {( 1 - z^{|v_l|})} {1-z} = $$

$$ 1 + \sum_{v_1 > ... > v_l > *} (-1)^l z^{|v_1| - |v_l| + 1}(1 + ... + z^{|v_l|-1}) = $$

$$ 1 + \sum_{v_1 > ... > v_l > *} (-1)^l (z^{|v_1| - |v_l| + 1} + ... + z^{|v_1|}) = $$

$$1 + \sum_{i\ge 0}(\sum_{a \ge i \ge a-b+1} s_{a,b})z^i.$$

\end{proof}

We may apply this result to compute $H(A(\Gamma),z)^{-1}$ 
where $\Gamma$ is the graph occuring in Cassidy and Shelton's example (Fig 1).  
The graph can be described as $\Gamma = (V,E)$ where $V = \cup_{i=0}^4 V_i$ with
$V_4 = \{u\},$ $ V_3 = \{w_1,w_2,w_3\},$ $ V_2 = \{x_1,x_2,x_3\},$ $ V_1 = \{y_1,y_2,y_3\},$ $ V_0 = \{*\}$ 
and there are edges from $u$ to $w_i$ for $i = 1,2,3$, 
from $w_i$ to $x_j$ if $i \ne j$, from $x_i$ to $y_j$ if $i \ne j$, 
and from $y_i$ to $*$ for $i = 1,2,3.$

\begin{figure}[htbp]\label{fig.pic}
\includegraphics[height=4in]{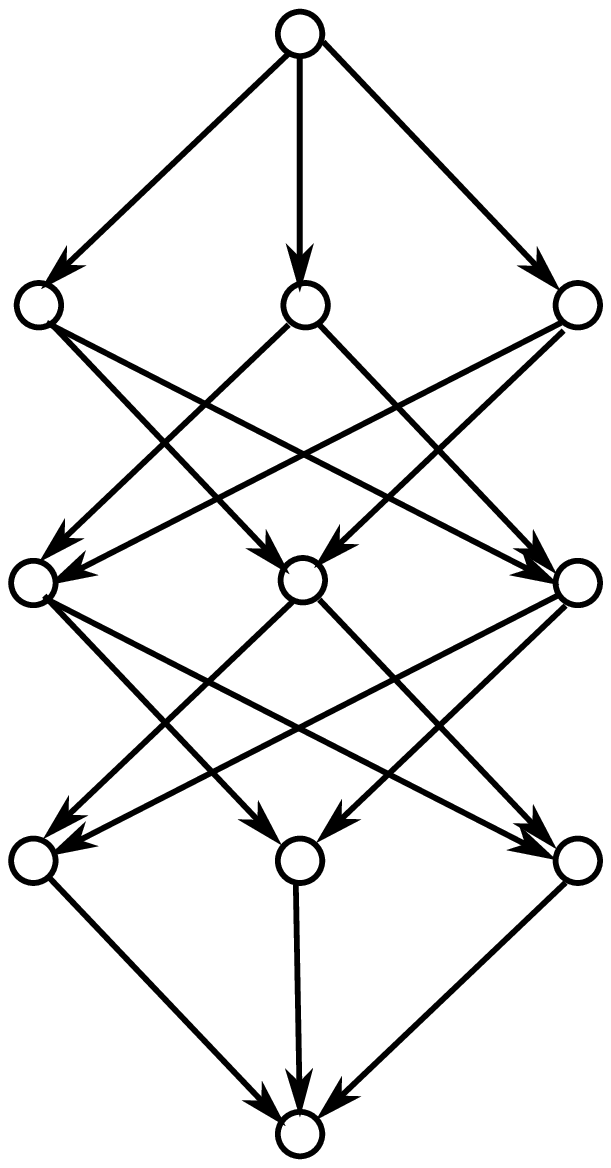}
\caption{ } 
\end{figure}

Then one sees that:
$$s_{4,4} = -1, s_{4,3} = 3, s_{4,2} = -3, s_{4,1} = 0,$$
$$s_{3,3} = -3, s_{3,2} = 6, s_{3,1} = -3,$$
$$s_{2,2} = -3, s_{2,1} = 6$$
and
$$s_{1,1} = -3.$$

Thus 
$$ H(A(\Gamma),z)^{-1} = - z^4 - z^3 + 8z^2 -10z +1.$$

Since $H(A(\Gamma),-z)^{-1}$ has a negative coefficient, 
it cannot be a Hilbert polyomial. Thus we recover 
Cassidy and Shelton's result (which they proved by homological methods) 
that  $A(\Gamma)$ is not Koszul.

\section{HILBERT SERIES OF $A(\Gamma)^!$}

\vskip 20 pt

We begin by noting the following description, which is immediate from 
the definition (cf \cite{SW}) of the Hilbert series of the 
dual $A^!$ of the quadratic algebra $A = T(V)/(<R_V>), R_V \subseteq V \otimes V.$ 
To simplify notation we will write $R$ in place of $R_V$. 
Write $V^k$ for the $k$-th tensor power of $V$ and define
$$R^{(0)} = K, R^{(1)} = V, R^{(2)} = R$$
and 
$$R^{(k+1)} = (V \otimes R^{(k)}) \cap (R \otimes V^{k-1})$$
for $k \ge 2.$

\begin{prop} If $A = T(V)/(<R>)$, then $A^!_j$ and $R^{(j)}$ are isomorphic vector spaces.
Consequently

$$H(A^!,z) = \sum_{j=0}^{\infty} dim(R^{(j)}) z^j.$$

\end{prop}

We will now investigate the $R^{(j)}$ that occur for the algebras $A(\Gamma).$

Let $\Gamma = (V,E)$ be a layered graph, $V = \sum_{i=0}^n V_i$.  Assume $v \in V_j, 1 \le j \le n$ and let
$l$ satisfy $1 \le l \le j.$  Denote by ${\mathcal C}(v,l)$ the set
$$\{v\otimes v_{(1)} \otimes ... \otimes v_{(l-1)}| v_{(m)} \in V_{j-m}, v_{(m)} > v_{(m+1)} \}  \subseteq T(V).$$
We call the elements of ${\mathcal C}(v,l)$ {\it linked monomials} of length $l$ starting at $v$.  
Also, if $w \in V_{j-l}$ define
$${\mathcal C}(v,w) = \{v\otimes v_{(1)} \otimes ... \otimes v_{(l-1)\otimes w} \in {\mathcal C}(v,l+1\}$$
and call elements of ${\mathcal C}(v,w)$ {\it linked monomials from $v$ to $w$}.

We say that a linked monomial $v\otimes v_{(1)} \otimes ... \otimes v_{(l-1)} \in {\mathcal C}(v,l)$ 
is an {\it admissible monomial} if, for each $p, 1 \le p \le l-1$ there is a linked monomial 
$$v\otimes w_{(1)} \otimes ... \otimes w_{(p-1)} \otimes  w_{(p)}\otimes w_{(p+1)}...\otimes w_{(l-1)} \in {\mathcal C}(v,l)$$
such that
$w_{(m)} = v_{(m)}$ for all $m, 1 \le m \le l-1, m \ne p$ 
and $v_{(p)} \ne w_{(p)}.$  
We denote the set of all admissible monomials in ${\mathcal C}(v,l)$ by ${\mathcal A}(v,l)$
and, if $w \in V_{j-l}$, define ${\mathcal A}(v,w) = {\mathcal A}(v,l+1) \cap {\mathcal C}(v,w).$

Clearly $R^{(k)} = V^{k-2}R \cap V^{k-3}RV \cap ... \cap RV^{k-2}$.  For $v \in V_j$ set
$R_v^{(k)} = vV^{k-1} \cap R^{(k)}.$

\begin{lem} $R_v^{(k)} \subset \ span \ {\mathcal A}(v,k).$

\end{lem}
\begin{proof} For $w \in V$ define a linear map
$$\phi_w:KV \rightarrow K$$
by 
$$\phi_w(u) = \delta_{w,u} \ \forall \ u \in V.$$

Let $$D = \{(w_1,w_2) \in V \times V|w_1 \not> w_2 \ or \ |w_1| - |w_2| \ne 1\}.$$
Then if $(w_1,w_2) \in D$ we have
$$\phi_{w_1} \otimes \phi_{w_2}R = 0.$$

Thus $$R^{(k)} \subset \bigcap_{l = 0}^{k-2} \bigcap_{(w_1,w_2) \in D} \ ker(id^{l}\otimes \phi_{w_1} \otimes \phi_{w_2} \otimes id^{k-2-l}).$$
Clearly 
$$vV^{k-1} \cap \bigcap_{l = 0}^{k-2} \bigcap_{(w_1,w_2 ) \in D} \ ker(id^{l}\otimes \phi_{w_1} \otimes \phi_{w_2}\otimes id^{k-2-l}) = \ span \ {\mathcal C}(v,k)$$

and so 
 
$$R_v^{(k)} \subset \ span \ {\mathcal C}(v,k).$$

Now define $\psi = \sum_{w \in V} \phi_w$ so that

$$\psi(u) = 1 \ \forall \ u \in V.$$

Then $R^{(2)}_v = {\mathcal C}(v,2) \cap \ker ( id \otimes \psi)$ and 
so

$$R^{(k)}_v =  {\mathcal C}(v,k) \cap (\bigcap_{l=1}^{k-l-1} \ ker (id^{l} \otimes \psi \otimes id^{k-l-1})).$$

But 
$$span \ {\mathcal C}(v,k) \cap \bigcap_{l=1}^{k-l-1} \ ker (id^{l} \otimes \psi \otimes id^{k-l-1}) \subseteq \  span \ {\mathcal A}(v,k)$$
and the lemma is proved.
\end{proof}

Assume $v \in V_j, w \in V_{j-l}, 2 \le l \le j$.  Define $R^{(l)}_{v,w} = R^{(l)} \cap {\mathcal A}(v,w).$

We say that the graph $\Gamma$ is ${\it oriented}$ if for any $v \in V_j, w \in V_{j-l}, 1 \le k \le l-1$ there exists a
function
$$o_{v,w}: {\mathcal A}(v,w) \rightarrow \{\pm1\}$$
such that if $m = v \otimes v_{(1)} \otimes ... \otimes v_{(l-1)}$ 
and $m' = v \otimes w_{(1)} \otimes ... \otimes w_{(l-1)} \in {\mathcal A}(v,w)$ satisfy
$v_q = w_q$ for $1 \le q \le l-1, q \ne k, and  v_k \ne w_k$ then
$$o(m) + o(m') = 0.$$
(This is similar to the definition of an oriented combinatorial cell complex given by Basak \cite{Basak}.)
Clearly there is at most one such $m'$.  However, $m \in {\mathcal A}(v,w)$ implies there that there is at least one such $m'$.  Thus $m'$ is unique.  We call it the {\it $k$-conjugate} of $m$.

Assume $\Gamma$ is oriented and  $v \in V_j, w \in V_{j-l}, 1 \le k \le l-1$.  Define the permutation
$$\tau_{v,w,k}: {\mathcal A}(v,w) \rightarrow {\mathcal A}(v,w)$$
to be the map that takes any $m \in {\mathcal A}(v,w)$ to its $k$-conjugate. Let ${\mathcal P}(v,w)$ denote the group of permutations of ${\mathcal A}(v,w)$ generated by
$\{\tau_{v,w,k}|1 \le k \le l-1\}.$
Clearly
$$o_{v,w}\tau_{v,w} = - \tau_{v,w}.$$
Extend elements of ${\mathcal P}(v,w)$ to endomorphisms of $span \ {\mathcal A}(v,w)$. 

\begin{lem} Let $\Gamma$ be oriented, $v, w \in V, v > w, |v| - |w| = l,$ and $1 \le k \le l-1.$  Let $u \in span \ {\mathcal A}(v,w)$.  
Then $(id^{k-1} \otimes \psi \otimes id^{l-k})u = 0 $ if and only if
$\tau_{v,w,k}u = -u.$
\end{lem}

\begin{proof}
For any linked monomial $m = y_{(1)} \otimes ... \otimes y_{(r)}$ define 
$$\psi_m = \phi_{y_{(1)}} \otimes ... \otimes \phi_{y_{(r)}}: V^{r} \rightarrow K.$$

For $m = v \otimes v_{(1)} \otimes ... \otimes v_{(j)} \otimes ... \otimes v_{(l-1)}\in {\mathcal C}(v,w), |v| - |w| = l$ define
$m^{[j]}_- = (id^j \otimes \psi^{l-j})m = v \otimes v_{(1)} \otimes ... \otimes v_{(j-1)}$, $m^{[j]} = v_{(j)}$, and
$m^{[j]}_+ = (\psi^{j+1} \otimes id^{l-j-1})m = v_{(j+1)} \otimes ... \otimes v_{(l-1)}.$
Thus $m = m^{[j]}_- \otimes m^{[j]} \otimes m^{[j]}_+.$

Let $m_1, m_2 \in {\mathcal A}(v,w)$.  Then
$m_{1,+}^{[k]} \otimes m_{1,-}^{[k]} = m_{2,+}^{[k]} \otimes m_{2,-}^{[k]}$ if and only if $m_1 = m_2$ or $m_1 = \tau_{v,w,k}m_2.$  Since the set of distinct monomals in $V^{l-1}$ is linearly independent,
$ker(id^k \otimes \psi \otimes id^{l-k-1}$ is spanned by
$\{m - \tau_{v,w,k}m|m \in {\mathcal A}(v,w)\}.$
Thus $ker(id^k \otimes \psi \otimes id^{l-k-1}) = ker(id + \tau_{v,w,k})$.

\end{proof}

Define $$m_{v,w} = \sum_{m \in {\mathcal A}(v,w)} o_{v,w}(m)m.$$

Let $X(v,w)$ denote a complete set of representatives for the orbits of 
${\mathcal A}(v,w)$ under ${\mathcal P}(v,w)$.
\begin{prop}
Assume that $\Gamma$ is oriented.  Then $\{m_{v,w}|m \in X(v,w)\}$ 
is a basis for ${\mathcal R}^{(l)}_{v,w}.$

\end{prop}
\begin{proof}
Since $$(id + \tau_{v,w,k})m_{v,w} = \sum_{m \in {\mathcal A}(v,w)} o_{v,w}(m) (m + \tau_{v,w,k}m) = $$
$$\sum_{m \in {\mathcal A}(v,w)}(o_{v,w}(m)+ o_{v,w} (\tau_{v,w,k}m)m =0,$$ 
the lemma shows that each $m_{v,w}$ is an element of ${\mathcal R}^{(l)}_{v,w}.$  
Moreover, the $m_{v,w}$ for $m \in X_{v,w}$ are sums over disjoint subsets of
${\mathcal A}(v,w)$ and so $\{m_{v,w}|m \in X(v,w)\}$ is linearly independent.

Now let $u = \sum_{m \in {\mathcal A}(v,w)} a(m)m \in {\mathcal R}^{(l)}_{v,w}.$  Then
$$u' = u - \sum_{m \in X(v,w)} a(m)m_{v,w} \in {\mathcal R}^{(l)}_{v,w}$$ 
and if we write
$u' = \sum_{m \in {\mathcal A}(v,w)}a'(m)m$, we have $a'(m) = 0 $ for all $m \in X(v,w).$  
Now by the lemma we see that $a'(m) = 0$ implies $a'(\tau_{v,w,k}m) = 0$ for all $k, 1 \le k \le l-1$.  Thus
$a'(\tau m) = 0$ for all $\tau \in {\mathcal P}(v,w)$ and $m \in X(v,w),$ 
so $u' = 0$ and the proposition is proved.

\end{proof}

\section{SMALL DIFFERENCE IN LEVELS}

Let $V$ be a finite-dimensional vector space over a field $K$.  Recall that
we write $V^k$ 
for the $k$-th tensor power
$V \otimes ... \otimes V \subseteq T(V)$.  
Let $R$ be a subspace of $V^2$, let $A = \sum_{n \ge 0} A_n$ 
be the (graded) quadratic algebra $T(V)/<R>$ and let 
$A^! = \sum_{n \ge 0} A^!_n$ be the quadratic dual of $A$.


\begin{prop}  Write $H(A,t)H(A^!,-t) = \sum_{n \ge 0} b_nt^n.$  
If $b_4 = 0$, then the lattice  in $V^4$ generated by $V^2R, VRV$ and $RV^2$ is distributive.
\end{prop}   
\begin{proof}
The lattice is distributive if and only if 
$V^2R\cap (VRV + RV^2) = V^2R\cap VRV + V^2R\cap RV^2.$
But $$V^2R\cap (VRV + RV^2) \supseteq V^2R\cap VRV + V^2R\cap RV^2$$ 
and so it is sufficient to prove that 
$$dim(V^2R\cap (VRV + RV^2)) = dim(V^2R\cap VRV + V^2R\cap  RV^2).$$

Write $R^{(j)} = V^{j-2}R \cap V^{j-3}RV \cap ... \cap RV^{j-2}$ 
and recall (Proposition 2.1) the well-known fact that
$dim(A^!_j) = dim \ R^{(j)}.$  

Now $$b_4 = dim(A_4) - dim(A_3)dim(A^!_1) + dim(A_2)dim(A^!_2) - $$
$$dim(A_1)dim(A^!_3) + dim(A^!_4) = $$
$$dim(V^4) - dim(V^2R + VRV + RV^2) - (dim(V^3) - dim(VR + RV))dim(V) + $$ 
$$(dim(V^2) - dim(R))dim(R) -
dim(V)dim(R^{(3)}) + dim(R^{(4)}).$$

Thus if $b_4 = 0$ we have
$$(dim(R))^2 + dim(V)dim(R^{(3)}) - dim(R^{(4)}) = $$
$$-dim(V^2R + VRV + RV^2) + dim(VR + RV))dim(V) + dim(V^2)dim(R).$$

Now the left-hand side may be written as 
$$dim(V^2R \cap RV^2) + dim(V^2R \cap VRV) - dim((V^2R \cap RV^2)\cap (V^2R \cap VRV))
=$$
$$dim((V^2R \cap RV^2) + (V^2R \cap VRV)).$$

Similarly, the right-hand side may be written as
$$- dim(V^2R + VRV + RV^2) + dim(VRV + RV^2) + dim(V^2R) =$$
$$ dim((V^2R)\cap(VRV + RV^2)).$$
Thus the proposition is proved.
\end{proof}

\section{CELL COMPLEXES}
Let $\mathcal C$ be a finite dimensional cell complex (see \cite{St}, Example 3.8.7) and
$M$ be the underlying manifold. Denote by
$P$ the ranked poset of cells of $\mathcal C$, ordered by defining $\sigma
_i\leq \sigma _j$ if $\bar \sigma _i\subseteq \bar \sigma _j$ with $rk
(\sigma )=dim (\sigma)+1.$ Let $\hat P$ be the ranked poset obtained
from $P$ by adding to $P$ the underlying space $M$ and the empty set
where $rk (M)=dim (M)+1$ and $rk(\emptyset)=0$. Recall (\cite{St}, Proposition
3.8.8) that $$ \mu (M, \emptyset)=\chi (M)-1$$
where $\mu (M,\emptyset )$ is the M\"obius function and $\chi (M)$ the
Euler characteristic of $M$. In the notations of \cite{St} $\mu (M, \emptyset)=\mu
_{\hat P}(\hat 0, \hat 1)$.

Let $\Gamma $ be the layered graph corresponding to the ranked poset $\hat
P$ and let $n = dim(M) + 2.$  Since the cells of $\mathcal C$ are 
connected, 
the layered graph $\Gamma$ is uniform. Therefore $A(\Gamma)$ is a quadratic algebra.  
Recall that $R$ denotes the space of defining relations.

\begin{thm} If the algebra $A(\Gamma )$ is numerically Koszul then
$$(-1)^ndim R^{(n)} = \chi (M)-1.$$
\end{thm}

\begin{proof} According to our computation \cite{RSW1, RW} (see also 
Lemma 1.3),
$$H(A(\Gamma),z)^{-1}=\mu (M, \emptyset)z^n +\dots +1.$$

Since the Hilbert polynomial $H(A(\Gamma)^!,z) = \sum_{j=0}^n dim(R^{(j)})z^j,$
if $ A(\Gamma)$ is numerically Koszul, then $\mu(M,\emptyset) = (-1)^ndim(R^{(n)})$ 
and the theorem follows from Proposition 3.8.8 of \cite{St}.

\end{proof}

\section{2-DIMENSIONAL MANIFOLDS}
\vskip 20 pt  
To any subdivision of a two-dimensional
connected, oriented manifold $M$ into faces, edges and vertices there corresponds a
layered graph  $\Gamma =(V,E)$, where $V_0=\{\emptyset \}$, $V_1$
is the set of vertices, $V_2$ is the set of edges, $V_3$ is the
set of faces, and $V_4=\{M\}$. An edge $x\in E$ goes from $a\in V_i$
to $b\in V_{i-1}$, $i=4,3,2,1$ if and only if $b$ belongs to the
border of $a$.

Set $g = |V_1|, h=|V_2|, f = |V_3|$ and $u=g + f$. Note that $u-h$ is the Euler
characteristic of $M$.
\medskip
\begin{prop} Suppose that in a subdivision of $M$ every
edge separates two distinct faces and the boundary of every edge consists
of two points. Then
$$H(A(\Gamma ),z)^{-1}=1- (1+u+h)z + (3h-1)z^2 - (h+1)z^3 +(u-h-1)z^4,$$
$$H(A(\Gamma )^!,z)=1 + (1+u+h)z + (3h-1)z^2 + (u-1)z^3 + z^4.$$
\end{prop}
The proof of this proposition will depend on a sequence of lemmas.

\begin{lem}  $\Gamma$ may be given the structure of an oriented graph.
\end{lem}

\begin{proof} For each edge $e$ we arbitrariy label the two endpoints as $h(e)$ 
and $t(e)$.  Let $F$ be any face of $M$ and $e$ be an edge of $F$.  
Then the orientation of $M$ induces an orientation of $F$ and 
this orientation of $F$ induces an orientation of the edge $e$.  
Then we define the orientation of the linked monomial 
$m_{F,e} = M \otimes F \otimes e$ by 
$o(m_{F,e}) = 1$ if the oriented edge $e$ is directed from $t(e)$ 
to $h(e)$ and $o(m_{F,e}) = -1$ if the oriented edge $e$ is directed 
from $h(e)$ to $t(e)$.  Let $v$ be an endpoint of $e$. We  
define the orientation of the linked monomial 
$m_{F,e,v} = F \otimes e \otimes v$ by $o(m_{F,e,v}) = 1$ if the 
oriented edge $e$ is directed towards $v$ and define $o(m_{F,e,v}) = -1$ 
if the oriented edge $e$ is directed away from $v$.  

Now for any edge $e$ there are exactly two faces incident to $e$ 
and the orientations induced on $e$ from $F_1 $ and $F_2$ are opposite.  
Also, if $F$ is any face and $v$ is any endpoint occuring in $F$, 
there are exactly two edges of $F$ incident to $v$ and one of 
these edges is directed towards $v$ and one is directed away from $v$.  
Thus the function $o$ we have defined gives $\Gamma$ the structure of an oriented graph.

\end{proof} 

\begin{lem}
(a)  Let $F$ be a face of $M$.  Then ${\mathcal P}(F,*)$ acts transitively on ${\mathcal A}(F,*).$ 

(b)  Let $w$ be a vertex of $M$.  Then ${\mathcal P}(M,w)$ acts transitively on ${\mathcal A}(M,w).$ 

(c)  ${\mathcal P}(M,*)$ acts transitively on ${\mathcal A}(M,*).$

\end{lem}

\begin{proof}

(a)  Let the face $F$ have vertices $v_1,...,v_l$ and edges $e_1,...,e_l$ 
with $e_i$ incident to $v_i$ and $v_{i+1}$ where the subscripts are taken modulo $l$. 
Then ${\mathcal A}(F,*) = \{F \otimes e_i \otimes v_i, F \otimes e_i \otimes v_{i+1}|1 \le i \le l\}, 
\tau_{F,*,1}$ interchanges $F \otimes e_i\otimes v_i$ and $F \otimes e_{i-1} \otimes v_i$ and 
$\tau_{F,*,2}$ interchanges $F \otimes e_i\otimes v_i$ and $F \otimes e_i \otimes v_{i+1}$.
Then $\tau_{F,*,2}\tau_{F,*,1}$ maps $F \otimes e_i \otimes v_i$ to $F \otimes e_{i-1} \otimes v_{i-1}$ and so the transitivity of ${\mathcal P}(F,*)$ acting on ${\mathcal A}(F,*)$ is clear. 

(b) Let the vertex $w$ be incident to edges $e_1,...,e_l$ and faces $F_1,...,F_l$ 
with $F_i$ incident to $e_i$ and $e_{i+1}$ where the subscripts are taken modulo $l$. 
Then ${\mathcal A}(M,w) = \{M \otimes F_i \otimes e_i, M \otimes F_i \otimes e_{i+1}|1 \le i \le l\}, 
\tau_{M,w,1}$ interchanges $M \otimes F_i\otimes e_i$ and $M \otimes F_{i-1} \otimes e_i$ and 
$\tau_{M,w,2}$ interchanges $M \otimes F_i\otimes e_i$ and $M \otimes F_i \otimes e_{i+1}$.
Then $\tau_{M,w,2}\tau_{M,w,1}$ maps $M \otimes F_i \otimes e_i$ to $M \otimes F_{i-1} \otimes e_{i-1}$ and so the transitivity of ${\mathcal P}(M,w)$ acting on ${\mathcal A}(M,w)$ is clear. 

(c) Let $M \otimes F_1 \otimes e_1 \otimes v_1$ and 
$M \otimes F_2 \otimes e_2 \otimes v_2$ be two elements of ${\mathcal A}(M,*)$.
Then since $M$ is connected, there is a path $\pi$ in $M$ from $v_1$ to $v_2$.  
We may assume that this path does not contain any endpoint other than $v_1$ and $v_2$.  
Let $n(\pi)$ denote the number of faces that $\pi$ traverses.

If $n(\pi ) = 1,$ then $F_1 = F_2$ and the result follows from the transitivity of ${\mathcal P}(F_1,*)$ on ${\mathcal A}(F_1,*).$

Now assume $n(\pi) > 1$.  Then the path $\pi$ remains in $F_1$ 
until it crosses some edge $e$ separating $F_1$ from some face $F'$.  
Let $v$ be an endpoint of $e$.  Then by the transitivity of   ${\mathcal P}(F_1,*)$ 
on ${\mathcal A}(F_1,*)$ there is some $\tau \in {\mathcal P}(M,*)$ 
such that $\tau(M \otimes F_1 \otimes e_1 \otimes v_1) = M \otimes F_1 \otimes e \otimes v.$  
But there is a path $\pi'$ from $v$ to $v_2$ with $n(\pi') = n(\pi) - 1.$  
By induction we may assume that there is 
$\tau' \in {\mathcal P}(M,*)$ such that 
$\tau'(M \otimes F_1 \otimes e \otimes v) = M \otimes F_2 \otimes e_2 \otimes v_2.$
Then $\tau'\tau( M \otimes F_1 \otimes e_1 \otimes v_1) = M \otimes F_2 \otimes e_2 \otimes v_2$ 
and the proof is complete.

\end{proof}

Let $G_{\Gamma}$ denote the set of all linear functionals $\alpha: KE \rightarrow K$ such that
for every face $F \in V_3$ we have
$$\sum o(F,e)\alpha(e) = 0$$
where the sum is taken over all edges $e$ incident to $F$.  
Let $F' \in V_3$ be a face with more than $3$ vertices and 
let $v,v'$ be vertices of $F'$ that are  not connected by an edge of $F'$.  
Then we may create a new subdivision of $M$ by dividing the face $F'$ 
into two faces $F''$ and $F'''$ by adding a new edge $e'$ connecting $v$ and $v'$.  
Let $\Gamma'$ be the graph corresponding to this new subdivision and 
let $E' = E \cup \{e'\}$, the set of edges of the new subdivision.  
Then any $\alpha \in G_{\Gamma}$ has a unique extension to an element
$\alpha' \in G_{\Gamma'}$ since we must have $\alpha'(e') = - \sum o(F,e)\alpha(e)$ 
where the sum is taken either over all edges $e$ of $F''$ different 
from $e'$ or, equivalently, over all edges $e$ of $F'''$ different from $e'$.  Define
$$\theta: (KE')^* \rightarrow (KE)^*$$
by
$$\theta: \beta \mapsto \beta|_{KE}.$$
Then $\theta(G_{\Gamma'}) = G_{\Gamma}.$

Let $w \in V_1$.  Define $q_{w,\Gamma} \in (KE)^*$ by
$$q_{w,\Gamma} = \sum_{t(e) = w} e^* - \sum_{h(e) = w} e^*.$$
Clearly $q_{w,\Gamma} \in G_{\Gamma}$ and $\theta(q_{w,\Gamma'}) = q_{w,\Gamma}$.

\begin{lem}
$G_{\Gamma} \subseteq \sum_{w \in V_1} Kq_{w,\Gamma}.$

\end{lem}
\begin{proof}
By the previous remarks, if
$$G_{\Gamma'} \subseteq \sum_{w \in V_1} Kq_{w,\Gamma'}$$
then
$$G_{\Gamma} \subseteq \sum_{w \in V_1} Kq_{w,\Gamma}.$$

Therefore, it is sufficient to prove the lemma under the additional assumption that every face is a triangle.

Let $\alpha \in G_{\Gamma}$ and suppose that $S \subseteq V_1$ has the property that 
$\alpha(e) = 0$ whenever both endpoints of $e$ are contained in $S$.  
Clearly $S = \emptyset$ has this property and if $S = V_1$ 
then $\alpha \in  \sum_{w \in V_1} Kq_{w,\Gamma}$.  Thus if, 
whenever $S \ne V_1$, we can find $S' \subseteq V_1$ properly 
containing $S$ and $\alpha' \in \alpha + \sum_{w \in V_1} Kq_{w,\Gamma}$ 
such that $\alpha'(e) = 0$ whenever both endpoints of $e$ 
are contained in $S'$, the lemma will follow by induction.

Now assume $S \ne V_1$. Since $M$ is connected, there is 
some $u \in V_1$ and some edge $e' \in V_2$ connecting $u$ to an element $s \in S.$
If $u = t(e')$ set $\alpha' = \alpha - \alpha(e')q_{u,\Gamma}$ and 
if $u = h(e')$ set $\alpha' = \alpha + \alpha(e')q_{u,\Gamma}$.
Then $\alpha'(e') = 0$ and $\alpha'(e) = 0$ for any edge $e$ connecting two vertices in $S$. Now,
if there is an edge $e''$ from $u$ to some vertex $s' \in S, s' \ne s$ then $u,s,s'$ 
are the vertices of some face and consequently $\alpha'(e'') = 0$.  Then we may take $S' = S \cup \{u\}$.

Next suppose there is no edge from $u$ to any element of $S$ other than $s$.  
Then, since $u$ is the endpoint of at least two edges, there is an edge $e'$ from $u$ to some $u' \notin S$.  
Then $u,u',s$ are the three vertices of some face $F'$ and so there is an edge $e''$ from $u'$ to $s$.  
If there is also an edge from $u'$ to some other element of $S$, then the previous argument 
with $u$ replaced by $u'$ shows that we may find the required $S'$.  
If there is no edge from $u'$ to a vertex of $S$ other than $s$, then setting $S' = S \cup \{u,u'\}$ gives the result.
\end{proof}

We can now prove the proposition.  

\begin{proof}
Because $V_4 = \{M\}$ and $V_0 = \{*\}$ we have
$$R^{(4)} = R^{(4)}_{M,*}.$$
Thus by Lemma 5.3 and Proposition 2.4  we have $$R^{(4)} = Km_{M,*}$$ and so $dim \ R^{(4)} = 1.$

Next, observe that $R^{(3)} = R^{(3)}_M \oplus \sum_{F \in V_3} R^{(3)}_{F,*}.$  
Again by Lemma 5.3 and Proposition 2.4, $R^{(3)}_{F,*} = Km_{F,*}$ 
and since $R^{(3)}_{F,*} \subseteq FV^2$ 
and $\sum_{F \in V_3} FV^2$ is direct, we have
$$dim(\sum_{F \in V_3} R^{(3)}_{F,*}) = |V_3| = f.$$

Now $dim(\sum_{w \in V_1} K(m_{M,w} \otimes w)) = |V_1|$. By Lemma 5.4
${ R}_M^{(3)} = \sum_{w \in V_1} {R}_{M,w}^{(3)}.$ 
Then,  by Lemma 5.3 and Proposition 2.4, the map
$$id^3 \otimes \psi: \sum K(m_{M,w} \otimes w) \rightarrow R^{(3)}_M$$
is surjective.  The kernel is
$$K{\mathcal A}(M,4) \cap (\cap_{i=1}^3 ker(id^i \otimes \psi \otimes id^{3-i}) = R^{(4)}_M.$$
Thus
$$0 \rightarrow R^{(4)}_M \rightarrow \sum K(m_{M,w}\otimes w) \rightarrow R^{(3)}_M \rightarrow 0$$
is exact and so $dim(R^{(3)}_M) = |V_1| - 1 = g-1.$  Thus 
$$dim(R^{(3)}) = f+g-1 = u-1.$$

Now $$R^{(2)} = R^{(2)}_M \oplus (\oplus _{F \in V_3} R_F^{(2)}) \oplus (\oplus _{e \in V_2} R^{(2)}_e).$$
Clearly
$dim  (R^{(2)}_M) = f-1.$
Since $dim \ R^{(2)}_F$ is equal to the number of edges of $F$ minus $1$ and since every edge is incident 
to two faces we have $dim(\oplus _{F \in V_3} R^{(2)}_F) = 2h-f.$  Finally, $dim 
(R^{(2)}_e) = 1$ and so
$dim(\oplus_{e \in V_2} R^{(2)}_e )= h.$  Thus $dim (R^{(2)}) = f-1 + 2h-f + h = 3h-1.$
As the number of generators of $A(\Gamma)^!$ is 
$|V_1| + |V_2| + |V_3| + |V_4| = g + h + f + 1 = u + h + 1$, the expression for
$H(A(\Gamma),z)^{-1}$ is proved.

The expression for $H(A(\Gamma)^!,z)$ follows from Lemma 1.3 by noting that
$$s_{4,4} = -1, s_{4,3} = f, s_{4,2} = -h, s_{4,1} = g,$$
$$s_{3,3} = -f, s_{3,2} = 2h, s_{3,1} = -2h,$$
$$s_{2,2} = -h, s_{2,1} = 2h,$$
and
$$s_{1,1} = -g.$$

\end{proof}

\medskip
\begin{cor} For $\Gamma $ as above, the following
statements are equivalent:

i) $A(\Gamma )$ is numerically Koszul;

ii) $A(\Gamma )$ is Koszul;

iii) the Euler characteristic of $M$ is $2$, i.e. $M$ is a topological
sphere.
\end{cor}
\begin{proof}
By Proposition 5.1, 
$H(A(\Gamma),Z)^{-1} - H(A(\Gamma)^!,-z) = (u-h-2)(z^3 + z^4).$  Thus (i) and (iii) are equivalent.  Also (i) implies (ii) by Proposition 3.1 and (iii) is well-known to imply (i).
\end{proof}

\section{KOSZULITY OF CERTAIN $A(\Gamma)$}

As noted in the introduction, we are indebted to Cassidy and Shelton 
for pointing out errors in our previous paper \cite{RSW1}.  First of all, Lemma 1.4 in that paper and the line immediately preceding it should be replaced by:  

For every $v \in V^+,$ set $\tilde v - v^{(0)} + ... + v^{(|v|-1)} \in FV^+.$ then 
$\tau(e_i) = e(v_{i-1},1) - e(v_i,1)$ and hence
$\nu(e_i) = \tilde v_{i-1} - \tilde v_i.$  Since 
$\theta$ is induced by $\nu$ we have

$$\theta(e(\pi,k)) = (-1)^k \sum_{1 \le i_1 < ... < i_k \le m} (\tilde v_{i_1 - 1} - \tilde v_{i_1})...(\tilde v_{i_{k-1}} - \tilde v_{i_k}).$$

A more serious error is that we have incorrectly asserted (in Theorems 4.6 and 5.2) 
that $gr \ A(\Gamma)$ and $A(\Gamma)$ are Koszul algebras 
for any uniform layered graph $\Gamma$ with a unique minimal vertex.  
The argument given for these assertions depends on a subsidiary result, 
Lemma 4.4, which asserts that if $l \ge 0, j \ge 1, v \in \cup_{j=2}^n V_j$ then
$$P_j(v)V^{l+1} \cap R^{(l+2)} = g_{l+3}(S_{j-1}(v)V^{l+2}\cap R^{(l+3)}).$$
Here $g_j:V^j \rightarrow V^{j-1}$ denotes $\psi \otimes id^{j-1}$ 
(recall that $\psi(v) = 1$ for all $v \in V$), $S_j(v)$ denotes $span\{w \in V|v > w, |v| - |w| = j\}$, 
and $P_j(v) = ker \ g_1 |_{S_j(v)}.$

The Cassidy-Shelton example provides a counter-example to this assertion.
Let $r = (x_1-x_2)(y_1-y_3) - (x_1-x_3)(y_1-y_2) = x_1(y_2-y_3) - x_2(y_1-y_3) + x_3(y_1-y_2)$, where 
the notation is as in Section 1. Then
$r \in  P_2(v)V \cap R^{(2)},$ while $S_1(u)V^3 \cap R^{(3)} \subseteq \sum_{i=1}^3 R_{w_i}^{(3)}  = (0).$

The proof in \cite{RSW1} shows: 
\begin{prop} If
$$P_j(v)V^{l+1} \cap R^{(l+2)} = g_{l+3}(S_{j-1}(v)V^{l+2}\cap R^{(l+3)})$$
whenever $l \ge 0, j \ge 1, v \in \cup_{j=2}^n V_j$ then
$gr \ A(\Gamma)$ and $A(\Gamma)$ are Koszul algebras.  
\end{prop}

We use this sufficient condition to show that  $A(\Gamma)$ is Koszul for certain classes of graphs $\Gamma$.

We begin by recalling (from \cite{RSW2}) that a {\it complete layered graph} 
is a layered graph $\Gamma = (V,E)$ with
$V = \cup_{i=0}^n V_i$ such 
that if $1 \le i \le n, v \in V_i, w \in V_{i-1}$ then there is an edge from $v$ to $w$.

\begin{thm} Let $V = (V,E)$ be a complete layered graph where $V = \cup_{i=0}^n V_i$ 
and $V_0 = \{*\}$.  Then $gr \ A(\Gamma)$ and $A(\Gamma)$ are Koszul algebras.

\end{thm}
\begin{proof} It is sufficient to prove the result for $gr \ A(\Gamma)$.  For 
simplicity of notation we will write $R^{(j)}$ for $(gr \ R)^{(j)}.$
Let $v \in V_k, j \le k$.  We first note that
$$vV^{j-1} \cap R^{(j)} = $$
$$span\{v(v_{k-1,1} - v_{k-1,2})...(v_{k-j+1,1} - v_{k-j+1,2})|v_{i,1},v_{i,2} \in V_i , 0 \le i \le k-1 \}.$$
To see this note that it is clear that $vV^{j-1} \cap R^{(j)}$ contains the right-hand side.
Also, any element $q$ of $vV^{j-1} \cap R^{(j)}$ is a sum of elements in
$vwV^{(j-2)} \cap vV^{(j-1)}$ where $w  \in V_{k-1}$.  
Equality  now follows by applying induction to $V^{(j-1)}$ and using $q \in RV^{j-2}.$ 

Using this equality we have
$$P_j(v)V^{l+1} \cap R^{(l+2)} = $$
$$P_j(v) span\{(v_{k-j,1}-v_{k-j,2})...(v_{k-j-l,1} - v_{k-j-l,2})|$$
$$v_{i,1},v_{i,2} \in V_i , k-j-l \le i \le k-j\} = $$
$$g_{l+3}(S_{j-1}(v)V^{l+2}\cap R^{(l+3)}),$$
proving the theorem.
\end{proof}

Next we will prove the Koszulity of splitting algebras associated to abstract simplicial complexes.
We recall that an abstract simplicial complex on a set $X$ is a subset 
$\Delta$ of the power set ${\mathcal P}(X)$ of $X$ such that 
$\{x\} \in \Delta$ for all $x \in X$ and if $S \in \Delta$ 
and $T \subseteq S$, then $T \in \Delta$. 

Let $X$ be a finite set and $\Delta$ be an abstract simplicial complex on $X$.  
Define a layered graph
$\Gamma_{\Delta} = (V_{\Delta}, E_{\Delta})$ by $V_{\Delta,i} = \{S \in \Delta| \ |S| = i\}$ 
and $E_{\Delta} = \{e_{S,s}|S \in \Delta, s \in S\}$
where $t(e_{S,s}) = S, h(e_{S,s}) = S\setminus \{s\}.$

Note that $\Delta_n = {\mathcal P}(\{1,...,n\})$ is an abstract simplicial complex in $\{1,...,n\}$.  
The algebra $ A(\Gamma_{\Delta_n})$ is the algebra studied in \cite{GRW}.  This algebra is denoted
by $Q_n$ there and we will continue to denote it by $Q_n$.

As above, it is sufficient to prove the result for $gr \ A(\Gamma)$ and 
we will write $R^{(j)}$ for $(gr \ R)^{(j)}.$
Set
$$R_{\Delta} = span\{S((S \setminus \{s\}) - (S \setminus \{t\}))|S \in \Delta, s,t \in S\},$$

$${\mathcal R}_{\Delta} = T(V)R_{\Delta}T(V),$$ 
and 
$${\mathcal R}_{\Delta}^{(k)} = \bigcap_{j=0}^{k-2} V_{\Delta}^jR_{\Delta}V_{\Delta}^{k-2-j}.$$
For $S \in \Delta, |S| \ge k$, define
$${\mathcal R}_{\Delta,S}^{(k)} = SV_{\Delta}^{k-1} \cap {\mathcal R}_{\Delta}^{(k)}.$$
Then, since ${\mathcal R}_{\Delta} = \oplus_{S \in \Delta} SV_{\Delta} \cap {\mathcal R}_{\Delta}$, we have
$${\mathcal R}_{\Delta}^{(k)} = \bigoplus_{S \in \Delta} {\mathcal R}_{\Delta,S}^{(k)}.$$

The following result is immediate.

\begin{lem}
Let $S = \{s_1,...,s_i\} \in \Delta$.  The map
$$j \mapsto s_j, 1 \le j \le i$$
extends to an injection 
$$T(V_{\Delta_i}) \rightarrow T(V_{\Delta})$$
and this map restricts to an isomorphism of vector spaces  
$${\mathcal R}_{\Delta_i,\{1,...,i\}} \rightarrow {\mathcal R}_{\Delta,S}.$$
\end{lem}

\begin{thm}
Let $\Delta$ be an abstract simplicial complex.  Then $A(\Gamma_{\Delta})$ is a Koszul algebra.
\end{thm}
\begin{proof}
Let $S = \{s_1,...,s_i\} \in \Delta$.  Let  
$$\phi:{\mathcal R}_{\Delta_i,\{1,...,i\}} \rightarrow {\mathcal R}_{\Delta,S}$$
be the isomorphism constructed in the lemma.

Now let $A \subseteq \{1,...,n\}$.  Then Proposition 6.4 of \cite{GGRSW} together with the well 
known identification of ${\mathcal R}_{\Delta}^{(k)}$ with 
the space of elements of degree $k$ in $A(\Gamma_{\Delta})_k$ (proposition 2.1) allows us to describe
a basis for $AV_{\Delta}^{k-1} \cap {\mathcal R}_{\Delta}^k.$  
For $B = \{b_1,...,b_k\} \subseteq A$ with
$b_1 < ... < b_k$ define
$$S(A:B) = \sum _{\sigma \in S_k} sgn(\sigma) A(A\setminus \{b_{\sigma(1)}\})...(A \setminus \{b_{\sigma(1)},...,b_{\sigma(k-1)}\}).$$
Then
$$\{\phi(S(A,B))|B \subseteq A \subseteq \{1,...,n\}, min \ A \notin B, |B| = k\}$$
is a basis for $AV_{\Delta}^{k-1} \cap {\mathcal R}_{\Delta}^k.$

Now if $B \subseteq A \subseteq \{1,...,n\}, min \ A \notin B, |B| = k$ then
$$g_kS(A:B) = \sum_{\sigma \in S_k} sgn(\sigma)(A \setminus \{b_{\sigma(1)}\})...(A \setminus \{b_{\sigma(1)},...,b_{\sigma(k-1)}\}).$$
Then writing $\sigma ' = \sigma \circ (12)$ and 
taking all sums over all $\sigma \in S_k$ such that $\sigma(1) < \sigma(2)$, we have
$$g_kS(A:B) = \sum sgn(\sigma)(A \setminus \{b_{\sigma(1)}\})...(A \setminus \{b_{\sigma(1)},...,b_{\sigma(k-1)}\}) + $$
$$\sum sgn(\sigma ')(A \setminus \{b_{\sigma'(1)}\})...(A \setminus \{b_{\sigma'(1)},...,b_{\sigma'(k-1)}\}) = $$
$$\sum sgn(\sigma)((A \setminus \{b_{\sigma(1)}\})- (A\setminus \{b_{\sigma(2)}\})
(A\setminus \{b_{\sigma(1)}, b_{\sigma(2)}\})...(A \setminus \{b_{\sigma(1)},...,b_{\sigma(k-1)}\}).$$
Thus $$g_kS(A:B) \in P_1(A)V^{k-2} \cap R^{(k)},$$
and so the hypothesis of Proposition 6.1 is satisfied and the theorem is proved.

\end{proof}

Since any Koszul algebra is numerically Koszul, we have the following corollary.
\begin{cor}  
Let $\Delta$ be an abstract simplicial complex.  Then $A(\Gamma_{\Delta})$ is a numerically Koszul algebra.
\end{cor}

We remark that this result may be established directly without the use of Theorem 6.4.

To see this, recall that $$H(A(\Gamma_{\Delta}),z) 
= (1-z)(1 + 
\sum_{v_1 > ... > v_l \ge *} (-1)^l z^{|v_1|-|v_l|+1})^{-1}$$

and that $$H(A(\Gamma_{\Delta})^!,z ) = \sum_{k \ge 0}  dim({\mathcal R}_{\Delta}^{(k)})z^k.$$
Thus $A(\Gamma_{\Delta})$ is numerically Koszul if and only if
$$(1-z)\sum_{k \ge 0} (-1)^k dim({\mathcal R}^{(k)})z^k = 1 + 
\sum_{v_1 > ... > v_l \ge *} (-1)^l z^{|v_1|-|v_l|+1}.$$

By Proposition 6.4 of \cite {GGRSW} we have that, if $i \ge k$, 
$$dim{\mathcal R}_{\Delta_i,\{1,...,i\}}^{(k)} = \binom {i-1} {k-1}.$$

Thus, in view of Lemma 6.1,  $$(1-z)\sum_{k \ge 0} (-1)^k dim({\mathcal R}^{(k)})z^k =
(1-z)\sum_{k \ge 0} (-1)^k \sum_{i \ge k} |V_{\Delta,i}|\binom {i-1} {k-1} z^k = $$
$$\sum_{k \ge 0} (-1)^k \sum_{i \ge k-1} |V_{\Delta,i}|\binom {i} {k-1} z^k.$$
  
Now consider 
$1 + 
\sum_{v_1 > ... > v_l \ge *} (-1)^l z^{|v_1|-|v_l|+1}.$  By Example 3.8.3 of \cite{St} we have 
$$\sum_{v_1 > ... > v_l \ge *} (-1)^l z^{|v_1|-|v_l|+1} = \sum_{v_1 > v_l \ge *}  (-z)^{|v_1|-|v_l|+1} =$$
$$\sum_{i \ge k} |V_{\Delta,i}|\binom {i} {k-1} (-z)^k,$$
establishing the numerical Koszulity of $A(\Gamma_{\Delta}).$

\begin{rem}
The layered graph $A(\Gamma)$ corresponding to an abstract simplicial complex will not, in general, have a unique maximal vertex.
Adding such vertex will destroy numerical Koszulity of the corresponding splitting algebra.
\end{rem}

\enddocument